\documentclass[]{article}
\usepackage[T1]{fontenc}
\usepackage[latin9]{inputenc}
\usepackage{color}
\usepackage{amstext}
\usepackage{amsmath,amsthm}
\usepackage{amssymb}
\usepackage{setspace}
\usepackage[a4paper,left=2.8cm,right=2.8cm,top=3cm,bottom=3.5cm]{geometry}

\usepackage[affil-it]{authblk}

\makeatletter
\theoremstyle{plain}	% title bold; body italic
\newtheorem{theorem}{Theorem}
\newtheorem{lemma}[theorem]{Lemma}
\newtheorem{coroll}[theorem]{Corollary}
\newtheorem{prop}[theorem]{Proposition}

\newtheorem*{claim*}{Claim}

\theoremstyle{definition}	% title bold; body not italic
\newtheorem{defn}[theorem]{Definition}

\newtheorem*{problem*}{Problem}

\theoremstyle{remark}
\newtheorem{remarkx}[theorem]{Remark}
\newenvironment{remark}
{\pushQED{\qed}\remarkx}
{\popQED\endremarkx}
\newtheorem{examplex}[theorem]{Example}
\newenvironment{example}
{\pushQED{\qed}\examplex}
{\popQED\endexamplex}
\newtheorem{fact}{Fact}

\usepackage{xcolor}
\usepackage{enumitem}
\newcommand{\defeq}{:=} %{\stackrel{\Delta}{=}}
\newcommand{\identity}{{\rm id}}
\newcommand{\dist}{{\rm dist}}
\newcommand{\gradient}{\,\mathrm{grad}\,}
\newcommand{\interior}{\mathrm{int\,}}
\newcommand{\cmnt}[1]{}
\newcommand{\red}[1]{{\leavevmode\color{red}#1}}

\usepackage{tikz}
\usetikzlibrary{calc,intersections,arrows.meta,bending,patterns,angles,matrix}
\usepackage{pgfplots}
\usepgfplotslibrary{fillbetween}

\usepackage{hyperref}
\hypersetup{
	colorlinks=false,
	allcolors=black
}

\title{On Wilson's theorem about domains of attraction and tubular neighborhoods}
\date{}
\author[]{Bohuan Lin\thanks{b.lin@rug.nl}}
\author[]{Weijia Yao\thanks{weijiayao1992@foxmail.com, corresponding author}}
\author[]{Ming Cao\thanks{m.cao@rug.nl}}
\affil[]{University of Groningen}

\begin{document}

%\begin{document}

\maketitle
%\begin{frontmatter}



%\title{On Wilson's theorem about domains of attraction and tubular neighborhoods}
%\author{Bohuan Lin\fnref{label1}}
%\ead{b.lin@rug.nl}
%\fntext[label1]{University of Groningen, the Netherlands}

% \affiliation{organization={},
%\author{Weijia Yao\corref{cor1}\fnref{label2}}
%\ead{weijiayao1992@foxmail.com}
% \fntext[label2]{National University of Defense Technology, China}

%\author{Ming Cao\fnref{label1}}
%\ead{m.cao@rug.nl}

\begin{abstract}
	In this paper, we show that the domain of attraction of a compact asymptotically stable submanifold of a finite-dimensional smooth manifold of an autonomous system is homeomorphic to its tubular neighborhood. The compactness of the attractor is crucial, without which this result is false; two counterexamples are provided to demonstrate this. %This result is a corrected version of Wilson's theorem \cite[Theorem 3.4]{wilson1967structure}.
	
\end{abstract}
\noindent\textbf{Keywords:\space} domain of attraction,  \space compact manifold, \space asymptotic stability, \space  autonomous system

%Research highlights
%\item We show a fundamental result about the domain of attraction of a compact asymptotically stable submanifold.
%\item We provide counterexamples to emphasize the key role of the compactness of the submanifold.
%\item The aforementioned result refines the existing one and  facilitates the topological analysis of autonomous systems.

%\begin{keywords}
%domain of attraction, \sep compact manifold, \sep asymptotic stability, \sep autonomous system



%\end{keywords}

%\end{frontmatter}

\section{Introduction}
The domain of attraction of an attractor of a continuous dynamical system has been widely studied. An \emph{attractor} is a closed invariant set of which there exists an open neighborhood such that every trajectory of the dynamical system starting within the neighborhood eventually converges to the attractor, in the sense that the distance between the trajectory and the attractor converges to zero; namely, the attractor is \emph{attractive}. And the set of all initial conditions rendering the corresponding trajectories to converge to the attractor is called the \emph{domain of attraction} of the attractor \cite{khalil2002nonlinear,bhatia2002stability}. Generally, it is difficult or sometimes impossible to find analytically the domain of attraction of an attractor. Since an attractor is attractive, if additionally it is \emph{Lyapunov stable} \cite[Chapter 4]{khalil2002nonlinear}, then it is called an \emph{asymptotically stable attractor}; sometimes Lyapunov functions can be utilized to estimate its domain of attraction, but the estimate can be conservative \cite[Chapers 4 and 8]{khalil2002nonlinear}. 

Partly due to the difficulty of calculating the domain of attraction of an attractor, some studies in the literature instead investigate the ``shapes'' or ``sizes'' of domains of attraction in the topological sense \cite{sontag2013mathematical,bernuau2019topological,moulay2010topological,bhatia2002stability,bhat2000topological,wilson1967structure}. In particular, in the simplest case where the attractor is an asymptotically stable equilibrium point, it has been shown in \cite[Theorem 21]{sontag2013mathematical} that the domain of attraction is contractible. This result characterizes the ``shape'' of the domain of attraction, and it also implies the ``size'' of the domain of attraction. Namely, it leads to the topological obstruction that if the state space of the system is not contractible, then an equilibrium point cannot be stabilized globally  \cite[Corollary 5.9.3]{sontag2013mathematical}. Another topological obstruction is shown in \cite{bhat2000topological}, which states that the domain of attraction of an asymptotically stable equilibrium point cannot be the whole state space (i.e., global asymptotic stability of an equilibrium is impossible) if the state space of the continuous dynamical system has the structure of a vector bundle over a compact manifold. Some studies partly generalize these results to asymptotically stable attractors that are not necessarily equilibrium points. In \cite{moulay2010topological}, it is proved that a compact, asymptotically stable attractor defined on a manifold (or more generally, on a locally compact metric space) is a \emph{weak deformation retract} of its domain of attraction. The conclusion is further developed in \cite{bernuau2019topological}, which shows that if the considered manifold is the Euclidean space $\mathbb{R}^n$, then the compact asymptotically stable attractor is a \emph{strong deformation retract} of its domain of attraction. 

Assuming that the asymptotically stable attractors are \emph{compact submanifolds} of some ambient finite-dimensional smooth manifolds, stronger conclusions can be made about the domains of attraction. For example, it is proved in  \cite[Chapter V, Lemma 3.2]{bhatia2002stability} that the intersection of an $\epsilon$-neighborhood of the attractor and some sublevel set of a corresponding Lyapunov function (of which the existence is automatically guaranteed \cite{wilson1969smoothing}) is a deformation retract of the domain of attraction of the attractor. This result is refined in \cite{moulay2010topological,yao2021topo}, which conclude that the attractor itself is a strong deformation retract of its domain of attraction. Therefore, the attractor and its domain of attraction are homotopy equivalent. This result has practical significance.  For example, it facilitates the analysis regarding the existence of singular points and the possibility of global convergence of trajectories to desired paths in the vector-field guided path-following problem ffor robotic control systems \cite{yao2021topo}. Note that the results discussed in this paragraph are strengthened for \cite{wilson1967structure} the case where the attractor is an embedded submanifold, where Theorem 3.4 in \cite{wilson1967structure} claims that the domain of attraction is diffeomorphic to a tubular neighborhood of the attractor, which can be either a compact or non-compact submanifold. However, in this paper, we will show that the compactness of the attractor is crucial, without which such a claim becomes false. In addition, the proof of Theorem 3.4 in \cite{wilson1967structure} is very brief, only indicating the method without giving sufficient detail. In this paper, we will detail the proof for a corrected version of this theorem, where the attractor is required to be compact. 

\textit{Contributions}: Throughout the paper, manifolds or submanifolds are \emph{without} boundaries, and they are second countable and paracompact. We assume that the attractor is compact, asymptotically stable and it is a submanifold of some finite-dimensional smooth manifold. We show that the compactness of the attractor is crucial for Theorem 3.4 in \cite{wilson1967structure} by providing counterexamples where Theorem 3.4 in \cite{wilson1967structure} no longer holds if the attractor is \emph{not} compact. Taking the compactness of the attractor into account, Theorem 3.4 in \cite{wilson1967structure} is corrected as below:
\begin{theorem}[Corrected version of Theorem 3.4 in \cite{wilson1967structure}]	\label{thm:DA homeomorphic to tub} 
	The domain of attraction of a compact asymptotically stable submanifold of a finite-dimensional smooth manifold	of an autonomous system is homeomorphic to its tubular neighborhood.
\end{theorem}
In this paper, we will give a complete and detailed proof of Theorem \ref{thm:DA homeomorphic to tub}, along with some auxiliary results to gain more insight into the theorem. 

The remainder of the paper is organized as follows. Section \ref{sec_pre} provides some preparatory results for the convenience of proving Theorem \ref{thm:DA homeomorphic to tub}. Then the detailed proof of Theorem \ref{thm:DA homeomorphic to tub} is elaborated in Section \ref{sec_proof}. To justify the importance of the compactness of the attractor in this theorem, we provide two counterexamples where the attractor is \emph{not} compact and hence Theorem \ref{thm:DA homeomorphic to tub}  fails to hold  in Section \ref{sec_example}. Finally, Section \ref{sec_conclu} concludes the paper.

\section{Preparatory results} \label{sec_pre}
In this section, we go through some basic notions and facts that will be used in the sequel. Let $\mathcal{M}$ and $\mathcal{N}$ be smooth manifolds, and $\mathcal{S}$ be a submanifold of $\mathcal{M}$. Note that in this section, the submanifold $\mathcal{S}$ can be compact or non-compact unless its compactness is specified explicitly. The notation $\defeq$ means ``defined to be''. The map $\identity$ is the identity map where the domain and codomain are clear from the context.

First, we recall the definitions of topological and smooth embeddings.
\begin{defn}[Topological and smooth embeddings, {\cite[p. 85]{lee2015introduction}}]
	A \emph{(topological) embedding} is an injective continuous map that is a homeomorphism onto its image (with the subspace topology). A \emph{smooth embedding} is a smooth immersion that is also a (topological) embedding.
\end{defn}
If $f: \mathcal{M} \to \mathcal{N}$ is an embedding, the image $f(\mathcal{M})$ can be regarded as a homeomorphic copy of $\mathcal{M}$ inside $\mathcal{\mathcal{N}}$. If $f: \mathcal{M} \to \mathcal{N}$ is a smooth embedding, then it is both a topological embedding and a smooth immersion.%, rather than being a (topological) embedding that happens to be smooth.

For each $p\in \mathcal{M}$, denote by $T_{p}\mathcal{M}$ and $T_{p}\mathcal{S}$ the tangent
spaces respectively of $\mathcal{M}$ and $\mathcal{S}$ at $p$, and by $T\mathcal{M}$ and $T\mathcal{S}$
the tangent bundles. Note that $T\mathcal{S}$ can be regarded as a subbundle of $T\mathcal{M}$ in a natural way. 

\begin{defn}[Normal bundle]
The normal bundle $\mathcal{N}_{\mathcal{S}}$ of $\mathcal{S}$ in $\mathcal{M}$ is the quotient bundle $T_{\mathcal{S}}\mathcal{M}\big/T\mathcal{S} \defeq \bigsqcup_{p\in \mathcal{S}} (T_{p}\mathcal{M} / T_p \mathcal{S})$, where $\bigsqcup$ denotes the disjoint union.%, where $T_{\mathcal{S}}\mathcal{M}=\bigsqcup_{p\in \mathcal{S}}T_{p}\mathcal{M}$. 
\end{defn}

\begin{fact}[{\cite[Sections 6.1 and 7.1]{mukherjee2016differential}}]
Let $g$ be any Riemannian metric on $\mathcal{M}$. For each $p\in \mathcal{M}$, let
$\mathcal{N}_{p}$ be the orthogonal complement of $T_{p}\mathcal{S}$ in $T_{p}\mathcal{M}$ with
respect to $g$. Then $\bigsqcup_{p\in \mathcal{S}}\mathcal{N}_{p}$ is a subbundle of
$T_{\mathcal{S}}\mathcal{M}$ and it is isomorphic to $T_{\mathcal{S}}\mathcal{M}\big/T\mathcal{S}$. This gives another way of defining the normal bundle of $\mathcal{S}$ in $\mathcal{M}$.%, or say, an embedding of $N_{\mathcal{S}}$ into $T_{\mathcal{S}}\mathcal{M}$.
\end{fact}
\begin{fact}[{\cite[Section 5.1]{mukherjee2016differential}}]
For any vector bundle $\mathcal{E}$ over $\mathcal{S}$, (the image of) the zero section of $\mathcal{E}$ can be canonically identified with $\mathcal{S}$ via 
\begin{align*}
	\iota_{\mathcal{S}}: \bar{0}_{\mathcal{S}} \subseteq \mathcal{E} &\to \mathcal{S}  \\
	0_{x} &\mapsto x
\end{align*}
where $\bar{0}_{\mathcal{S}} \subseteq \mathcal{E}$ denotes (the image of) the zero section of $\mathcal{E}$, and $0_{x}$ denotes the zero vector in the vector space $\mathcal{E}_{x}$ for $x \in \mathcal{S}$. Therefore, $\iota_{\mathcal{S}}$ is a diffeomorphism from $\bar{0}_{\mathcal{S}}$ to $\mathcal{S}$. Note that viewing $\mathcal{S}$ as a submanifold of $\mathcal{M}$, $\iota_{\mathcal{S}}$ can also be regarded as an embedding of $\bar{0}_{\mathcal{S}}$ into $\mathcal{M}$.
\end{fact}
\begin{defn}[Tubular neighborhood] \label{def:tubular-neighborhood}
A tubular neighborhood of $\mathcal{S}$ is an open embedding $\tau:\mathcal{E} \rightarrow \mathcal{M}$ from some vector bundle $\mathcal{E}$ over $\mathcal{S}$ to $\mathcal{M}$ satisfying 
\[
	\tau\big|_{\bar{0}_{\mathcal{S}}}=\iota_{\mathcal{S}}.
\]
More loosely, we often call the open set $\mathcal{W} \defeq \tau(\mathcal{E})$ a tubular neighborhood of $\mathcal{S}$.% An open set $W$ is also called a tubular neighborhood of $\mathcal{S}$ if there exists such an embedding that $W=\tau(E)$.
\end{defn}
Whether we refer to a tubular neighborhood as an embedding or an open set should be clear from the context. 

\begin{theorem}[Existence of tubular neighborhood, {\cite[Proposition 7.1.3]{mukherjee2016differential}}] \label{thm:existence of tubular neighborhoods-1}
Suppose that $\mathcal{S}$ is a submanifold of $\mathcal{M}$. Then there exists an embedding
$\tau:N_{\mathcal{S}} \to \mathcal{M}$ from the normal bundle $N_{\mathcal{S}}$ of $\mathcal{S}$
into $\mathcal{M}$ such that $\tau$ keeps the zero section of $N_{\mathcal{S}}$ (i.e., $\tau(0_{x})=x$ for all $x\in\mathcal{S}$, or $\tau\big|_{\bar{0}_{\mathcal{S}}}=\iota_{\mathcal{S}}$).
\end{theorem}
\begin{remark}
	This means that $\tau: N_\mathcal{S} \to \mathcal{M}$ is a tubular neighborhood of $\mathcal{S}$, and $\tau$ is a diffeomorphism between $N_\mathcal{S}$ and $\tau(N_\mathcal{S})$.
\end{remark}

Before presenting the uniqueness result of tubular neighborhoods, we first recall the definitions of \emph{isotopy} and \emph{diffeotopy}.
\begin{defn}[Isotopy and diffeotopy, {\cite[pp. 177-178]{hirsch2012differential}}]
	An \emph{isotopy} from $\mathcal{M}$ to $\mathcal{N}$ is a map $F : \mathcal{M} \times \mathcal{I} \to \mathcal{N} $, where $\mathcal{I} \subseteq \mathbb{R}$ is an interval, such that for each $t \in \mathcal{I}$, the map
	$
	F_t: \mathcal{M} \to \mathcal{N}
	$
	defined by $ x \mapsto F(x,t)$ is an embedding. We also say $F$ is an isotopy from $F_0$ to $F_1$, and $F_0$ and $F_1$ are called \emph{isotopic}.
	If each $F_t$ is a smooth embedding, then $F$ is a \emph{smooth isotopy} from $\mathcal{M}$ to $\mathcal{N}$.
	If each $F_t$ is a diffeomorphism, \cmnt{\red{and $F_0=\identity_\mathcal{N}$}, where $\identity_\mathcal{N}$ denotes the identity map of $\mathcal{N}$,} then $F$ is called a \emph{diffeotopy}.
\end{defn}
Throughout the paper, whenever we mention an \emph{isotoopy}, we mean a \emph{smooth isotopy}. Now we show the uniqueness result of the tubular neighborhood as follows.

\begin{theorem}[Uniqueness of tubular neighborhood I, {\cite[Theorem 7.4.4]{mukherjee2016differential}}] \label{thm:uniqueness of tubular neighborhoods-1}
Suppose that $f_{i}:\mathcal{E}_{i}\rightarrow\mathcal{M}, \;i=0,1$, are tubular neighborhoods of $\mathcal{S}$. Then there exists a bundle map\footnote{More specifically, the bundle map $\gamma$ is a bundle isomorphism. This is because $f_0$ and $f_1$ are embeddings and their images are open sets in $\mathcal{M}$; therefore, $\mathcal{E}_0$ and $\mathcal{E}_1$ are vector bundles which, as manifolds,  have the same dimensions as $\mathcal{M}$ does. } $\lambda:\mathcal{E}_{0}\rightarrow \mathcal{E}_{1}$
such that $f_{0}$ and $f_{1}\circ\lambda$ are isotopic (see  Fig. \ref{fig: thm:uniqueness of tubular neighborhoods-1}).
\end{theorem}
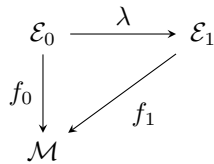
\begin{figure}[h]
	\centering
	\begin{tikzpicture}
	\matrix (m) [matrix of math nodes,row sep=3em,column sep=4em,minimum width=2em]
	{
		{\mathcal{E}_0} & {\mathcal{E}_1} \\
		{\mathcal{M}} & {} \\};
	\path[-stealth]
	(m-1-1) edge node [above] {$\lambda$} (m-1-2) 
	(m-1-1) edge node [left] {$f_0$} (m-2-1) 
	(m-1-2) edge node [below right] {$f_1$} (m-2-1);
	\end{tikzpicture}
	\caption{Relations in Theorem \ref{thm:uniqueness of tubular neighborhoods-1}. If $f_0$ and $f_1$ are two tubular neighborhoods, then there exists a bundle map $\lambda: \mathcal{E}_0 \to \mathcal{E}_1$ such that $f_{0}$ and $f_{1}\circ\lambda$ are isotopic.}
	\label{fig: thm:uniqueness of tubular neighborhoods-1}
\end{figure}
Due to Theorems \ref{thm:existence of tubular neighborhoods-1} and \ref{thm:uniqueness of tubular neighborhoods-1}, Theorem \ref{thm:DA homeomorphic to tub} implies the following result. 
\begin{prop}
	The domain of attraction of a compact asymptotically stable submanifold $\mathcal{S}$ of a finite-dimensional smooth manifold $\mathcal{M}$ of an autonomous system is homeomorphic to its normal bundle $N_{\mathcal{S}}$.
\end{prop}
\begin{proof}
	Combine Theorems  \ref{thm:DA homeomorphic to tub}, \ref{thm:existence of tubular neighborhoods-1} and \ref{thm:uniqueness of tubular neighborhoods-1}.
\end{proof}

Denote by $G: \mathcal{E}_0 \times (-\delta, 1+\delta) \to \mathcal{M}$ the isotopy from $f_{0}$ to $f_{1}\circ\lambda$. Then Theorem \ref{thm:uniqueness of tubular neighborhoods-1} implies that $G_t: \mathcal{E}_0 \to \mathcal{M}$ is a tubular neighborhood for any $t \in (-\delta, 1+\delta)$. Now let 
	\[
		h(x,t)=G(f_{0}^{-1}(x),t)
	\]
	for $(x,t)\in f_0(\mathcal{E}_0) \times(-\delta,1+\delta)$. We have the following corollary. 
\begin{coroll}[Uniqueness of tubular neighborhood II, {\cite[Theorem 7.4.4]{mukherjee2016differential}}] \label{cor:uniqueness of tubular neighborhoods}  
Suppose that $\mathcal{S}$ is a submanifold of $\mathcal{M}$, and $\mathcal{W}_{0}$ and $\mathcal{W}_{1}$ are two tubular neighborhoods (as open sets) of $\mathcal{S}$ in $\mathcal{M}$, then there exists an isotopy $h:\mathcal{W}_{0}\times(-\delta,1+\delta)\rightarrow\mathcal{M}$ such that 
\[
	h_{0}=j_{\mathcal{W}_{0}}, \quad h_{1}(\mathcal{W}_{0})=\mathcal{W}_{1}, \quad h_{t}\big|_{\mathcal{S}}=j_{\mathcal{S}}
\]
for every $t\in(-\delta,1+\delta)$, where $h_t \defeq h(\cdot, t)$, $j_{\mathcal{W}_{0}}$ and $j_{\mathcal{S}}$ are the inclusions of $\mathcal{W}_{0}$ and $\mathcal{S}$ into $\mathcal{M}$
respectively.
\end{coroll}
%\begin{remark}
	Therefore, any two tubular neighborhoods $\mathcal{W}_0$ and $\mathcal{W}_1$ are homeomorphic.
%\end{remark}

\begin{defn}[Closed tubular neighborhood]
Fix a Euclidean metric $g$ on the vector bundle $\mathcal{E}$ over $\mathcal{S}$, and for any $r>0$, let\footnote{Note that $\mathcal{BE}_{r}$ is a submanifold of $\mathcal{E}$ with boundary $\partial(\mathcal{BE}_{r})=\{v\in \mathcal{E} : g(v,v)=r^{2}\}$.}
\[
	\mathcal{BE}_{r}=\{v\in \mathcal{E} : g(v,v)\leq r^{2}\}.
\]
A \emph{closed tubular neighborhood} $\mathcal{K}$ of $\mathcal{S}$ is a closed neighborhood of $\mathcal{S}$ in $\mathcal{M}$ such that there is an embedding $\phi:\mathcal{BE}_{r}\rightarrow \mathcal{M}$ satisfying
\[
	\phi(\mathcal{BE}_{r}) = \mathcal{K}, \quad \phi\big|_{\bar{0}_{\mathcal{S}}}=\iota_{\mathcal{S}}.
\]
\end{defn}

\begin{remark} \label{remark: BErEr}
If $\mathcal{S}$ is compact, then $\mathcal{BE}_{r}$ is by definition a closed tubular neighborhood of $\bar{0}_{\mathcal{S}}$ in $\mathcal{E}$ and that it is compact. Since $\mathcal{E}_{r}=\{v\in \mathcal{E} : g(v,v)<r^{2}\}$ can homeomorphically map to $\mathcal{E}$ while keeping the zero section, it is an (open) tubular neighborhood of $\bar{0}_{\mathcal{S}}$ in $\mathcal{E}$. In particular, $\mathcal{E}$ itself is a tubular neighborhood of $\bar{0}_{\mathcal{S}}$ in $\mathcal{E}$.
\end{remark}

Due to Remark \ref{remark: BErEr}, the following proposition holds.
\begin{prop} \label{prop:precompact tubular neighborhoods}
If $\mathcal{S}$ is compact, then there exists some tubular neighborhood $\mathcal{W}$ of $\mathcal{S}$ such that its closure $\bar{\mathcal{W}}$ is a closed tubular neighborhood which is also compact.
\end{prop}

We will use a technique which relies on the following results to prove Theorem \ref{thm:DA homeomorphic to tub} later.
\begin{lemma}[{\cite[Chapter 8, Theorem 1.4]{hirsch2012differential}}] \label{thm: isotopyextension} 
Suppose that $\mathcal{U}$ is an open set of the manifold $\mathcal{N}$ and that $\mathcal{S}$ is a compact subset of $\mathcal{N}$ contained in $\mathcal{U}$. Suppose that $h: \mathcal{U} \times(-\delta,1+\delta)\rightarrow\mathcal{N}$ is an isotopy with $h_{0}:\mathcal{U} \rightarrow \mathcal{N}$ being the inclusion. Then for any $\delta'\in(0,\delta)$, there exists a diffeotopy $H:\mathcal{N} \times(-\delta',1+\delta')\rightarrow\mathcal{N}$ with some open neighborhood $\mathcal{U}_{0}$ of $\mathcal{S}$ in $\mathcal{U}$ such that
\[
	H\big|_{\mathcal{U}_{0}\times(-\delta',1+\delta')}=h\big|_{\mathcal{U}_{0}\times(-\delta',1+\delta')}.
\]
\end{lemma}
% \begin{proof}
% 	See Appendix \ref{app1}.
% \end{proof}
%
\begin{remark}
Let $\tilde{h}$ be the level preserving map\footnote{The map $\tilde{h}$ is called the \emph{track} of $h$ \cite[p. 111]{hirsch2012differential}.}:
\begin{align*}
	\tilde{h}: \mathcal{U} \times(-\delta,1+\delta) &\to \mathcal{N}\times(-\delta,1+\delta) \\
	(p,t) &\mapsto \big(h_{t}(p),t\big).
\end{align*}
Note that Theorem 1.4 in Chapter 8 of \cite{hirsch2012differential} requires $\tilde{h}\big( \mathcal{U} \times(-\delta,1+\delta)\big)$ to be open in $\mathcal{N}\times(-\delta,1+\delta)$. However, this requirement is unnecessary at least in our case, since it can be easily checked that $\tilde{h}$ is a submersion\footnote{This is because $\tilde{h}$ is an immersion and the dimensions of $\mathcal{U}$ and $\mathcal{N}$ are the same.} and hence an open map.
\end{remark}

\begin{coroll} \label{cor:isotopyextension}
Suppose that $\mathcal{U}$ is an open set of the manifold $\mathcal{N}$ and that $\mathcal{S}$ is a compact subset of $\mathcal{N}$ contained in $\mathcal{U}$. Suppose that $h':\mathcal{U}\times(-\delta,1+\delta)\rightarrow\mathcal{N}$
is an isotopy, and there exits a diffeomorphism $f_{0}:\mathcal{N} \to \mathcal{N}$ that agrees with $h'_{0}$ on $\mathcal{U}$; i.e., 
\begin{equation} \label{eq_f0u}
	f_0|_{\mathcal{U}} = h_0'.
\end{equation}
Then for any $\delta'\in(0,\delta)$, there is a diffeotopy $F: \mathcal{N} \times(-\delta',1+\delta')\rightarrow\mathcal{N}$ with some open neighborhood $\mathcal{U}_{0}$ of $\mathcal{S}$ in $\mathcal{U}$ such that 
\[
	F\big|_{\mathcal{U}_{0}\times(-\delta',1+\delta')}=h'\big|_{\mathcal{U}_{0}\times(-\delta',1+\delta')}.
\]
\end{coroll}
\begin{proof}
Let $h=f_{0}^{-1}\circ h'$. Therefore, from \eqref{eq_f0u}, we have $h_0=f_{0}^{-1}\circ h'_0=j_{\mathcal{U}}$, where $j_{\mathcal{U}}: \mathcal{U} \to \mathcal{N}$ is the inclusion map from $\mathcal{U}$ to $\mathcal{N}$.  According to Lemma \ref{thm: isotopyextension}, there is a diffeotopy $H$ such that $H\big|_{\mathcal{U}_{0}\times(-\delta',1+\delta')}=h\big|_{\mathcal{U}_{0}\times(-\delta',1+\delta')}$. Then let $F=f_{0}\circ H$.
\end{proof}
\begin{remark}
Note that the open set $\mathcal{U}$ in the theorems above may be $\mathcal{N}$ itself, which is the case in Lemma \ref{lem:(extension-of-tubular-neighborhoods)} to be discussed later.
\end{remark}

Now we prove a lemma concerning tubular neighborhoods of the submanifold $\mathcal{S}$ of $\mathcal{M}$. This lemma greatly facilitates the arguments in Section \ref{sec_proof}.

Note that $(N_{\mathcal{S}}, \pi, \bar{0}_{\mathcal{S}})$, where $\pi: N_{\mathcal{S}}\rightarrow\bar{0}_{\mathcal{S}}$ defined by $p \mapsto 0_{x}$ for any $p\in N_{x}$ and $x\in \mathcal{S}$, is a vector bundle over $\bar{0}_{\mathcal{S}}$. Though this might be trivial since $\bar{0}_{\mathcal{S}}$ is identical to $\mathcal{S}$ in a canonical way, we still point it out as follows for the sake of clarity from the set-theoretic perspective.
\begin{figure}[h]
	\centering
	\begin{tikzpicture}
	\matrix (m) [matrix of math nodes,row sep=3em,column sep=4em,minimum width=2em]
	{
		{N_\mathcal{S}} & {N_\mathcal{S}} \\
		{N_\mathcal{S}} & {N_\mathcal{S}} \\};
	\path[-stealth]
	(m-1-1) edge node [right] {$j$} (m-2-1) 
	(m-1-1) edge node [above] {$\lambda$} (m-1-2) 
	(m-1-2) edge node [right] {$\identity_{N_\mathcal{S}}$} (m-2-2);
	\end{tikzpicture}
	\caption{Proof of Lemma \ref{lem:(extension-of-tubular-neighborhoods)}. }
	\label{fig: lem:(extension-of-tubular-neighborhoods)}
\end{figure}
\begin{lemma}[Extension of tubular neighborhoods] \label{lem:(extension-of-tubular-neighborhoods)} 
Suppose that $j:N_{\mathcal{S}}\rightarrow N_{\mathcal{S}}$ is a tubular
neighborhood of $\bar{0}_{\mathcal{S}}$; i.e., $j$ is an embedding and
$j(0_{x})=0_{x}$ for all $x\in \mathcal{S}$. Then for any compact set $\mathcal{K}$
in $N_{\mathcal{S}}$, there is a diffeomorphism $\beta$ on $N_{\mathcal{S}}$ such that
$\beta$ agrees with $j$ on some neighborhood of $\mathcal{K}$.
\end{lemma}

\begin{proof}
The idea is to use Corollary \ref{cor:isotopyextension}. To this end, we seek an isotopy $h:N_{\mathcal{S}}\times(-\delta,1+\delta)\rightarrow N_{\mathcal{S}}$ such that $h_{1}=j$ and $h_{0}$ is a diffeomorphism on $N_{\mathcal{S}}$. 

Note that both $\identity_{N_{\mathcal{S}}}:N_{\mathcal{S}}\rightarrow N_{\mathcal{S}}$ and $j:N_{\mathcal{S}}\rightarrow N_{\mathcal{S}}$
are tubular neighborhoods of $\bar{0}_{\mathcal{S}}$ in $N_{\mathcal{S}}$. Hence, according
to Theorem \ref{thm:uniqueness of tubular neighborhoods-1}, there exists a bundle isomorphism $\lambda:N_{\mathcal{S}}\rightarrow N_{\mathcal{S}}$ such
that there exists an isotopy $h$ from $\identity_{N_{\mathcal{S}}}\circ\lambda$
to $j$ (see Fig. \ref{fig: lem:(extension-of-tubular-neighborhoods)}). Since $\identity_{N_{\mathcal{S}}}\circ\lambda$ is a diffeomorphism, according
to Corollary \ref{cor:isotopyextension}, there exists a diffeotopy
$H:N_{\mathcal{S}}\times(-\delta,1+\delta)\rightarrow N_{\mathcal{S}}$ such that $H$
agrees with $h$ on some neighborhood of $\mathcal{K}$. Let $\beta=H_{1}$
and then it is a diffeomorphism and agrees with $h_{1}=j$ on such
a neighborhood.
\end{proof}

\section{Proof of Theorem \ref{thm:DA homeomorphic to tub}} \label{sec_proof}
The proof of Theorem \ref{thm:DA homeomorphic to tub} is based on \cite[Lemma 3.3]{wilson1967structure}. For clarity, we decompose the proof into several propositions. Denote by $\mathcal{M}$ the state space with a vector field $X$. Denote by
	$\varphi$ the flow of $X$ and assume that $\mathcal{S}$ is a \emph{compact} boundaryless
	submanifold of $\mathcal{M}$ and is an asymptotic
	stable attractor of $\varphi$. Denote by $\mathcal{D}_{A}$ the domain of attraction of
	$\mathcal{S}$. 

We start by fixing a precompact tubular neighborhood 
\[
	f_{o}:N_{\mathcal{S}}\rightarrow\mathcal{W}
\]
of $\mathcal{S}$ in $\mathcal{D}_{A}$, where $\mathcal{W} \defeq f_o(N_{\mathcal{S}})$. The existence of $f_o$ is guaranteed by Proposition \ref{prop:precompact tubular neighborhoods}.
\begin{prop} \label{prop:compactsets in tubularneighborhoodbyflow}
For each compact set $\mathcal{K}$ in the domain of attraction $\mathcal{D}_{A}$ , there exists
some $T_{K}>0$, such that $\varphi^{T}(\mathcal{K})\subseteq\mathcal{W}$ for
any $T>T_{K}$. Consequently, $\mathcal{K} \subseteq \varphi^{-T}(\mathcal{W})$ for any $T>T_{K}$.
\end{prop}

\begin{proof}
Due to the asymptotic stability of $\mathcal{S}$, there is some neighborhood $\mathcal{U}$ of $\mathcal{S}$ in $\mathcal{W}$ such that $\varphi^{[0,\infty)}(\mathcal{U})\subseteq\mathcal{W}$.
For any $x\in \mathcal{K}$, there is some $T_{x}>0$ with some neighborhood $\mathcal{B}_{x}$ of $x$ such that $\varphi^{T_{x}}(\mathcal{B}_{x})\subseteq \mathcal{U}$. Since $\mathcal{K}$ is compact, there is $\{\mathcal{B}_{x_{i}} \} _{i=1,...,k}$, where $k < \infty$, such that $\bigcup_{i} \mathcal{B}_{x_{i}}\supseteq \mathcal{K}$. Let $T_{K} \defeq \max_{i=1,...,k} T_{x_{i}}$ and the proof is completed.
\end{proof}
Note that $\mathcal{S}$ is invariant under $\varphi$, and hence $\mathcal{S} \subseteq\varphi^{-T}(\mathcal{W})$
for any $T \in \mathbb{R}$. Since $\varphi^{-T}:\mathcal{W}\rightarrow \mathcal{W}_{T} \defeq \varphi^{-T}(\mathcal{W})$
is a diffeomorphism and $\mathcal{W}$ is a tubular neighborhood
of $\mathcal{S}$,  it is natural to conjecture that $\mathcal{W}_{T}$ should also be a tubular neighborhood
of $\mathcal{S}$. This is indeed true as shown in the next proposition, but it is not straightforward. According to Definition \ref{def:tubular-neighborhood}, we still need to find a diffeomorphism $f_{T}$ from $N_{\mathcal{S}}$ to $\mathcal{W}_{T}$
such that $f_{T}\big|_{\bar{0}_{\mathcal{S}}}=\iota_{\mathcal{S}}$. Although $f=\varphi^{-T}\circ f_{o}$
is a diffeomorphism from $N_{\mathcal{S}}$ to $\mathcal{W}_{T}$, we have $f\big|_{\bar{0}_{\mathcal{S}}}=\varphi^{-T}\circ\iota_{\mathcal{S}}$, which is not necessarily equal to $\iota_{\mathcal{S}}$, and hence $f:N_{\mathcal{S}}\rightarrow \mathcal{W}_{T}$ is not necessarily a tubular neighborhood.
Yet $f\big|_{\bar{0}_{\mathcal{S}}}=\varphi^{-T}\circ\iota_{\mathcal{S}}$ and $\iota_{\mathcal{S}}$
are isotopic as maps from $\bar{0}_{\mathcal{S}}$ to $\mathcal{W}_{T}$ while $f$ and
$\varphi^{T}$ are both diffeomorphisms. This makes it possible to
use Lemma \ref{thm: isotopyextension}.
\begin{prop}
\label{prop:tubular neighborhoods by the flow}For any $T>0$, $\mathcal{W}_{T} \defeq \varphi^{-T}(\mathcal{W})$
is a tubular neighborhood of $\mathcal{S}$ in $\mathcal{D}_{A}$. That is,
there exists a diffeomorphism $f_{T}:N_{\mathcal{S}}\rightarrow \mathcal{W}_{T}$ such
that $f_{T}\big|_{\bar{0}_{\mathcal{S}}}=\iota_{\mathcal{S}}$.
\end{prop}

\begin{proof}
Obviously $f=\varphi^{-T}\circ f_{o}$ is a diffeomorphism from $N_{\mathcal{S}}$
to $\mathcal{W}_{T}$ with $0_{x}\in\bar{0}_{\mathcal{S}} \mapsto \varphi^{-T}(x)$.
Now we need to ``rectify'' the map. Denote by $f_{\mathcal{S}}$ the restriction
of $f$ on $\bar{0}_{\mathcal{S}}$ . Then $j_{1}=f^{-1}\circ\varphi^{T}\circ f_{\mathcal{S}}$
is a map mapping $\bar{0}_{\mathcal{S}}$ diffeomorphically to $\bar{0}_{\mathcal{S}}$. Let $j_{s}=f^{-1}\circ\varphi^{s\cdot T}\circ f_{\mathcal{S}}$
for $s\in(-\delta,1+\delta)$ and then $j:\bar{0}_{\mathcal{S}}\times(-\delta,1+\delta)\rightarrow N_{\mathcal{S}}$
is an isotopy such that $j_{0}$ is the inclusion map, and $f\circ j_{1}=\iota_{\mathcal{S}}$
on $\bar{0}_{\mathcal{S}}$.

Note that $g=\varphi\circ f$ with $g(x,t)=\varphi^{t}\circ f(x)$
is a smooth map from $N_{\mathcal{S}}\times\mathbb{R}$ to $\mathcal{D}_{A}$.
Since $\varphi^{[-\delta,1+\delta]\cdot T}\circ f(\bar{0}_{\mathcal{S}})=\mathcal{S}\subseteq \mathcal{W}_{T}$
and $[-\delta,1+\delta]\cdot T$ is compact, there exists an open
neighborhood $\mathcal{U}$ of $\bar{0}_{\mathcal{S}}$ in $N_{\mathcal{S}}$ such that $\varphi^{[-\delta,1+\delta]\cdot T}\circ f(\mathcal{U})\subseteq \mathcal{W}_{T}$.
Moreover, for any fixed $s\in[-\delta,1+\delta]$, $\varphi^{s\cdot T}\circ f(\cdot)$
is an injective submersion, and hence a smooth embedding. Define 
\begin{align*}
	h&:\mathcal{U} \times(-\delta,1+\delta)\rightarrow N_{\mathcal{S}} \\
	h(x,s)&=f^{-1}\circ\varphi^{s\cdot T}\circ f(x),
\end{align*}
which is an isotopy with $h_{0}$ being the inclusion map of $\mathcal{U}$ into $N_{\mathcal{S}}$
and $h_{s}\big|_{\bar{0}_{\mathcal{S}}}=j_{s}$. Then by Lemma \ref{thm: isotopyextension}, there exists a diffeotopy $H:N_{\mathcal{S}}\times(-\delta',1+\delta')\rightarrow N_{\mathcal{S}}$ for $\delta'\in(0,\delta)$ such that $H$ agrees with $h$ on $\mathcal{U}_{0}\times(-\delta',1+\delta')$
for some open neighborhood $\mathcal{U}_{0}$ of $\mathcal{S}$. 

Let $f_{T}=f\circ H_{1}$ and this is a diffeomorphism between $N_{\mathcal{S}}$
and $\mathcal{W}_{T}$. Moreover, restricted on $\mathcal{S}$, $f_{T}=f\circ h_{1}=f\circ j_{1}=\iota_{\mathcal{S}}$. Hence, $f_{T}:N_{\mathcal{S}}\rightarrow \mathcal{W}_{T}$ is a tubular neighborhood.
\end{proof}
Since the domain of attraction $\mathcal{D}_{A}$ is a smooth manifold
with the second countability, there exists an ascending chain of compact subsets $\mathcal{K}_{0}\subseteq \mathcal{K}_{1}\subseteq \cdots$ such that $\bigcup_{i\in\mathbb{N}}\mathcal{K}_{i}=\mathcal{D}_{A}$.
Choose $0<T_{0}<T_{1}<\cdots$ such that 
\[
	\mathcal{W}_{i} \defeq \varphi^{-T_{i}}(\mathcal{W})
\]
contains $\mathcal{K}_{i}$ for each $i$ and that $\overline{\mathcal{W}}_{i}\subseteq \mathcal{W}_{i+1}$.
This is possible due to the precompactness of $\mathcal{W}$. By Proposition
\ref{prop:tubular neighborhoods by the flow}, there exist tubular
neighborhoods $f_{i}:N_{\mathcal{S}}\rightarrow \mathcal{W}_{i}$ for all $i\in\mathbb{N}$.
The strategy to prove Theorem \ref{thm:DA homeomorphic to tub}
is to construct by induction an ascending chain of compact subsets $\mathcal{C}_{0}\subseteq \mathcal{C}_{1}\subseteq \cdots$ with tubular neighborhoods
$g_{i}:N_{\mathcal{S}}\rightarrow \mathcal{W}_{i}$ ``rectified'' from $f_{i}$ such that
$g_{i}(\mathcal{C}_{i})\supseteq \mathcal{K}_{i}$, $g_{i+1}$ agrees with $g_{i}$ on $\mathcal{C}_{i}$
and $\bigcup_{i}\mathcal{C}_{i}=N_{\mathcal{S}}$. Then the theorem follows by defining
a map $g:N_{\mathcal{S}}\rightarrow\mathcal{D}_{A}$ with $g=g_{i}$ on $\mathcal{C}_{i}$.
\begin{theorem}
There exists a diffeomorphism $g:N_{\mathcal{S}}\rightarrow\mathcal{D}_{A}$
such that $g\big|_{\bar{0}_{\mathcal{S}}}=\iota_{\mathcal{S}}$. 
\end{theorem}
\begin{proof}
Let $\mathcal{K}_{0}\subseteq \mathcal{K}_{1}\subseteq \cdots$ be an ascending chain of compact subsets such that $\bigcup_{i\in\mathbb{N}}\mathcal{K}_{i}=\mathcal{D}_{A}$ and $\mathcal{K}_{0}\supseteq \mathcal{S}$. Since $\mathcal{W}$ is precompact in $\mathcal{D}_{A}$, $\varphi^{-T}(\mathcal{W})$ is precompact for any $T>0$ in $\mathcal{D}_{A}$. Then by Proposition \ref{prop:compactsets in tubularneighborhoodbyflow}
we can choose inductively $0<T_{0}<T_{1}< \cdots$ such that $\overline{\mathcal{W}}_{i}\cup \mathcal{K}_{i+1}\subseteq \mathcal{W}_{i+1}$.
According to Proposition \ref{prop:tubular neighborhoods by the flow},
for each $i\in\mathbb{N}$, there is a diffeomorphism $f_{i}:N_{\mathcal{S}}\rightarrow \mathcal{W}_{i}$
such that $f_{i}(0_{x})=x$ for all $x\in \mathcal{S}$. Now we construct $\{(g_{i},\mathcal{C}_{i}) : i\in\mathbb{N}\}$
with $\mathcal{C}_{i}$ being compact sets in $N_{\mathcal{S}}$ and $g_{i}:N_{\mathcal{S}}\rightarrow \mathcal{W}_{i}$
being tubular neighborhoods such that 
\begin{enumerate}[label = (\roman*)]
	\item\label{req1}  $\mathcal{C}_{i}\subseteq \interior \mathcal{C}_{i+1}$;
	\item\label{req2}  $g_{i}(\mathcal{C}_{i})\supseteq \mathcal{K}_{i}$;
	\item\label{req3}  $g_{i+1}\big|_{\mathcal{C}_{i}}=g_{i}\big|_{\mathcal{C}_{i}}$;
	\item\label{req4}  $\bigcup_{i\in\mathbb{N}}\mathcal{C}_{i}=N_{\mathcal{S}}$;
\end{enumerate}

Take $g_{0}=f_{0}$ and $\mathcal{C}_{0}=\mathcal{BE}_{r_{0}}$ with $r_{0}$ large enough such that $\mathcal{BE}_{r_{0}}\supseteq g_{0}^{-1}(\mathcal{K}_{0})$. Let $j_{1}=f_{1}^{-1}\circ g_{0}$.
Then $j_{1}:N_{\mathcal{S}}\rightarrow N_{\mathcal{S}}$ is a tubular neighborhood of
$\bar{0}_{\mathcal{S}}$ in $N_{\mathcal{S}}$ and $f_{1}\circ j_{1}=g_{0}$. According
to Lemma \ref{lem:(extension-of-tubular-neighborhoods)}, there is
a bundle isomorphism $\beta_{1}:N_{\mathcal{S}}\rightarrow N_{\mathcal{S}}$ such that
$\beta_{1}$ agrees with $j_{1}$ on $\mathcal{C}_{0}$. Let $g_{1}=f_{1}\circ\beta_{1}$.
Then $g_{1}:N_{\mathcal{S}}\rightarrow \mathcal{W}_{1}$ is a diffeomorphism and $g_{1}=g_{0}$
on $\mathcal{C}_{0}$. Take $r_{1}$ large enough such that $r_{1}>2r_{0}$
and $\mathcal{C}_{1}=\mathcal{BE}_{r_{1}}$contains $g_{1}^{-1}(\mathcal{K}_{1})$. 

Suppose that for $n\in\mathbb{N}$, $\mathcal{A}_{n}=\{(g_{i},\mathcal{C}_{i}) : 0\leq i\leq n\}$
such that \ref{req1}, \ref{req2}, \ref{req3} are satisfied and $\mathcal{C}_{n}=\mathcal{BE}_{r_{n}}$ with $r_{n}>2^{n}r_{0}$.
Let $j_{n+1}=f_{n+1}^{-1}\circ g_{n}$. Then $j_{n+1}:N_{\mathcal{S}}\rightarrow N_{\mathcal{S}}$
is a tubular neighborhood of $\bar{0}_{\mathcal{S}}$ in $N_{\mathcal{S}}$. Again, according
to Lemma \ref{lem:(extension-of-tubular-neighborhoods)}, there exists
a diffeomorphism $\beta_{n+1}$ on $N_{\mathcal{S}}$ such that $\beta_{n+1}=j_{n+1}$
on $\mathcal{C}_{n}$. Set $g_{n+1}=f_{n+1}\circ\beta_{n+1}$ and then $g_{n+1}=g_{n}$
on $\mathcal{C}_{n}$. Pick a positive number $r_{n+1}$ such that $r_{n+1}>2r_{n}$
and $\mathcal{C}_{n+1}=\mathcal{BE}_{r_{n+1}}\supseteq g_{n+1}^{-1}(\mathcal{K}_{n+1})$. Then $\mathcal{A}_{n+1}=\mathcal{A}_{n}\cup\{(g_{n+1},\mathcal{C}_{n+1})\}$
again satisfies \ref{req1}, \ref{req2}, \ref{req3} with $r_{n+1}>2^{n+1}r_{0}$. By induction
we have $\{(g_{k},\mathcal{C}_{k}) : k\in\mathbb{N}\}$ satisfying \ref{req1}, \ref{req2},
\ref{req3} with $r_{k}>2^{k}r_{0}$ for all $k\in\mathbb{N}$. 

Define $g:N_{\mathcal{S}}\rightarrow\mathcal{D}_{A}$ with $g\big|_{\mathcal{C}_{i}}=g_{i}\big|_{\mathcal{C}_{i}}$
for all $i\in\mathbb{N}$. Then $g$ is well defined and $\textrm{Im } g=\mathcal{D}_{A}$
due to \ref{req3} and \ref{req2} respectively. Moreover, since $g_{i}$'s are diffeomorphisms
and $\bigcup_{i} \interior \mathcal{C}_{i}=N_{\mathcal{S}}$, $g$ is a local diffeomorphism. It is also obvious that \ref{req4} is satisfied. For any $p,q\in N_{\mathcal{S}}$, there exists $i$ such that $\mathcal{C}_{i}$ contains
$p$ and $q$. Then $g(p)=g(q)$ $\implies$ $g_{i}(p)=g_{i}(q)$
$\implies$ $p=q$. Hence $g$ is also injective. Therefore, $g$
is a diffeomorphism from $N_{\mathcal{S}}$ onto $\mathcal{D}_{A}$. Moreover,
since $\mathcal{K}_{0}$ is chosen to contain $\mathcal{S}$ in the beginning and $g_{0}$
keeps the zero section (i.e., $g(0_{x})=x$ for all $x\in\mathcal{S}$), $\mathcal{C}_{0}\supseteq\bar{0}_{\mathcal{S}}$. Therefore $g\big|_{\bar{0}_{\mathcal{S}}}=g_{0}\big|_{\bar{0}_{\mathcal{S}}}=\iota_{\mathcal{S}}$,
which concludes the proof.
\end{proof}

\section{Two counterexamples} \label{sec_example}
In this section, we illustrate two counterexamples to invalidate the original claim in \cite[Theorem 3.4]{wilson1967structure}. The first counterexample in Section \ref{sec:count1} discusses how \cite[Theorem 3.4]{wilson1967structure} fails to hold when the attractor is a noncompact manifold. The idea of constructing the counterexample is straightforward, but it usually involves an \emph{incomplete} Riemannian manifold as the ambient space. Nevertheless, another counterexample in Section \ref{sec:count2} involves a \emph{complete} Riemmannian manifold as the ambient space. The idea of the counterexample is to present two topologically equivalent dynamical systems, where the domains of attraction of the noncompact attractors are not homotopy equivalent. As a result, the domain of attraction of the noncompact attractor of either of the system is of a different homotopy type from its tubular neighbourhood, contradicting \cite[Theorem 3.4]{wilson1967structure}. Note that all the vector fields of the dynamical systems in this section are complete; i.e., solutions exist for all $t \in \mathbb{R}$.

\subsection{$\mathcal{M}$ is an incomplete Riemannian manifold} \label{sec:count1}
Theorem 3.4 in \cite{wilson1967structure} states that the domain of attraction of a uniformly asymptotically stable attractor, be it a compact or non-compact manifold, of a complete autonomous system is diffeomorphic to its tubular neighborhood. While the argument in Section \ref{sec_proof} holds for a compact attractor $\mathcal{S}$, it does not hold for a \emph{non-compact} attractor, since Proposition \ref{prop:compactsets in tubularneighborhoodbyflow} may be invalid when the attractor is noncompact. More specifically, when $\mathcal{S}$ is noncompact, it is possible that none of its tubular neighborhood contains any $\epsilon$-neighborhood of $\mathcal{S}$. To see this, note that if we take out one point from a submanifold, the $\epsilon$-neighbourhood of the new submanifold will only miss one point compared to that of the original submanifold, while its tubular neighbourhood (viewed as a vector bundle) would lose the whole fiber over the missing point. Exploiting this observation, we can construct a counterexample by starting with a compact asymptotically stable attractor and then taking one fixed point out of it.
\begin{example}
    Start with the smooth function $\bar{f}(x)=(\dist(x,\mathbb{S}^1))^2$ on $\mathbb{R}^2$ and let $\bar{X}=-\gradient \bar{f}$. This system has the unit circle $\mathbb{S}^1 \subseteq \mathbb{R}^2$ as the asymptotically stable attractor, and all points on $\mathbb{S}^1$ are fixed points.
    Now consider the state space $\mathcal{M}=\mathbb{R}^{2}-\{(1,0)\}$. Let $\mathcal{S}=\mathbb{S}^1-\{(1,0)\}$. It is a closed set and also a submanifold of $\mathcal{M}$, but it is noncompact. Let $f$ be the function on $\mathcal{M}$ such that $f(x)=(\dist(x,\mathcal{S}))^2$. The function $f$ is the restriction of $\bar{f}$ on $\mathcal{M}$, and hence it is smooth. The vector field $X=-\gradient f$ is then the restriction of $\bar{X}$ on $\mathcal{M}$, and it has $\mathcal{S}$ as an attractor, which is uniformly asymptotically stable. The domain of attraction is $\mathcal{M}-\{(0,0)\}$, which is not contractible. However, a tubular neighborhood of $\mathcal{S}$ is homeomorphic to $\mathcal{S} \times \mathbb{R}$, which is contractible, contradicting Theorem 3.4 in \cite{wilson1967structure}.
\end{example}

\subsection{$\mathcal{M}$ is a complete Riemannian manifold} \label{sec:count2}

In this section we demonstrate a dynamical system $(\mathcal{M},\varphi)$ as
a counterexample to Theorem 3.4 in \cite{wilson1967structure}
where the state space $\mathcal{M}$ is a complete Riemannian manifold and the
asymptotically stable attraction $\mathcal{S}$ is \emph{not} compact. Instead of directly showing the construction of the flow map $\varphi$ on $\mathcal{M}$, we first construct an auxiliary
system $(\mathcal{M}_{0},\varphi_{0})$, and then obtain $(\mathcal{M},\varphi)$ via
a topological conjugacy \cite[Chapter 2]{brin2002introduction} $h:\mathcal{M}_{0}\rightarrow \mathcal{M}$ . As an extra benefit
to be seen later, such a demonstration shows that uniformly
asymptotic stability is rather a ``geometric'' concept than a ``topological'' one. Namely, even if two dynamical systems are topologically conjugate, properties concerning the uniform asymptotic stability of the systems may not be (fully) preserved by the conjugacy.

\subsubsection{The auxiliary system $(\mathcal{M}_{0},\varphi_{0})$}

Let
\[
\mathcal{M}_{0}=\{(x,y,z)\in\mathbb{R}^{3} : x^{2}+z^{2}=1\}
\]
and
\[
\mathcal{S}_{0}=\{(x,y,z)\in \mathcal{M}_{0} : x=0,z=1\}.
\]
Endow $\mathcal{M}_{0}$ with the Riemannian metric $g_{0}$ induced by the
standard Riemannian metric $(dx)^{2}+(dy)^{2}+(dz)^{2}$ on $\mathbb{R}^{3}$.
Then $(\mathcal{M}_{0},g_{0})$ is a complete Riemannian manifold with the distance
$d_{\mathcal{M}_{0}}$. 

Let $Y_{0},Z_{0}$ be the vector fields on $\mathcal{M}_{0}$ defined by
\[
Y_{0}(x,y,z)=\begin{cases}
e^{-\frac{1}{y}}\frac{\partial}{\partial y}\big|_{(x,y,z)} & y>0\\
0 & y\leq0
\end{cases}
\]
and 
\[
Z_{0}(x,y,z)=x\cdot \left( x\frac{\partial}{\partial z}-z\frac{\partial}{\partial x} \right) \Big|_{\mathcal{M}_{0}}.
\]

Let $X_{0}=Y_{0}+Z_{0}$ and denote by $\varphi_{0}$ the
flow of $X_{0}$ on $\mathcal{M}_{0}$. Then $\mathcal{S}_{0}$ is a uniformly
asymptotically stable manifold of the dynamical system $(\mathcal{M}_{0},\varphi_{0})$
with its domain of attraction being 
\[
\mathcal{D}_{0}=\{(x,y,z)\in \mathcal{M}_{0} : z>-1\}.
\]
The following characterization of the stability of $\mathcal{S}_{0}$ will be
needed later. Namely, given any $a>-1$ with $\mathcal{W}'_{z>a}=\mathcal{M}_{0}\cap\{(x,y,z)\in \mathcal{M}_{0} : z>a\}$, corresponding to each $\epsilon>0$, there exists some $T_{\epsilon}>0$
such that $d_{\mathcal{M}_{0}}\big(\varphi_{0}^{[T_{\epsilon},+\infty)}(\mathcal{W}'_{z>a}),\mathcal{S}_{0}\big)<\epsilon$.
To see this, denote by $(x'_{t},y'_{t},z'_{t})$ the orbit $\varphi_{0}^{t}(p')$
for $p'=(x',y',z')\in \mathcal{M}_{0}$. Then $(x'_{t},z'_{t})\subseteq \mathbb{S}^{1}$
is subject to the equation 
\begin{equation} \label{eq:x'(t)-z'(t)}
\frac{d}{dt}(x'_{t},z'_{t})=\big(-x'_{t}z'_{t}, {x'_{t}}^{2}\big).
\end{equation}
Note that the dynamical system (\ref{eq:x'(t)-z'(t)}) on $\mathbb{S}^{1}$
has the point $q'_{0}=(0,1)$ as an asymptotically stable equilibrium
with the domain of attraction $\{(x,z)\in \mathbb{S}^{1} : z\neq-1\}$. Hence
for any $\epsilon>0$, there exists $T'>0$ such that for any $t\geq T'$
and $q'=(x',z')\in \mathbb{S}^{1}$ with $z'\geq a$, $\dist\big(\phi^{t}(q'),q'_{0}\big)<\epsilon$,
where $\dist$ is the distance on $\mathbb{S}^{1}$ measured by lengths of minor arcs, and $\phi$ is the flow of \eqref{eq:x'(t)-z'(t)}. Therefore,
\[
d_{\mathcal{M}_{0}}\big(\varphi_{0}^{t}(x',y',z'),\mathcal{S}_0 \big)\leq \dist\big(\phi^{t}(x',z'),q'_{0}\big)<\epsilon
\]
for all $t>T'$ and $(x',y',z')\in W'_{z>a}$.

For a point $(0,y,-1)\in \mathcal{M}_{0}-\mathcal{D}_{0}$ with $y\leq0$, it holds that $X_{0}\big|_{(x,y,z)}=0$.
For any $y>0$,
\begin{equation}
\label{eq:X0-(0,y,-1)}
X_{0}\big|_{(0,y,-1)}=Y_{0}\big|_{(0,y,-1)}=e^{-\frac{1}{y}}\frac{\partial}{\partial y}\Big|_{(0,y,-1)},
\end{equation}
implying
\begin{equation}
\varphi_{0}^{t}(0,y,-1)=(0,\gamma(t),-1)
\end{equation}
with 
\begin{equation}
\label{eq:dot_gamma>0}
\dot{\gamma}(t)=e^{-\frac{1}{\gamma(t)}}>0.
\end{equation}
Therefore, both $\gamma(t)$ and $\dot{\gamma}(t)$ increase strictly with respect to $t>0$.

\subsubsection{The system $(\mathcal{M},\varphi)$}

Now we construct the dynamical system $(\mathcal{M},\varphi)$ which will serve
as a counterexample. More specifically, a vector field $X$ on some Riemannian manifold $(\mathcal{M},g)$ is to be constructed with a uniformly asymptotically stable submanifold $\mathcal{S}$ of which the domain of attraction $\mathcal{D}$ is not homotopy equivalent to $\mathcal{S}$ itself. 

Let 
\[
r(y)=\begin{cases}
1-e^{-\frac{1}{y}} & y>0\\
1 & y\leq0
\end{cases}.
\]
Let $\mathcal{M}$ be the two-dimensional cylinder embedded in $\mathbb{R}^{3}$
defined by 
\[
	\mathcal{M} = \{(x,y,z) \in \mathbb{R}^3 : x^{2}+z^{2}=r(y) \},
\]
and let 
\[
    \mathcal{S}=\{(x,y,z) \in \mathcal{M} : x=0,z=\sqrt{r(y)}\}.
\]
Then $\mathcal{S}$ is an embedded submanifold and a closed subset in $\mathcal{M}$. Endowed
with the Riemannian metric $g_{\mathcal{M}}$ induced by the standard Riemannian
metric $g=(dx)^{2}+(dy)^{2}+(dz)^{2}$ on $\mathbb{R}^{3}$, $\mathcal{M}$
is a complete Riemannian manifold. Note that although the Riemannian metric $g_{\mathcal{M}}$ is induced by $g$, the corresponding distance $d_{\mathcal{M}}$
on $\mathcal{M}$ is not the restriction on $\mathcal{M}$ of the Euclidean distance $d$
on $\mathbb{R}^{3}$. Generally speaking, it holds that $d_{\mathcal{M}}(p,q)\geq d(p,q)$
for $p,q\in \mathcal{M}$. However, the topology $\tau_{\mathcal{M}}$ induced by $d_{\mathcal{M}}$
on $\mathcal{M}$ is exactly the subspace topology inherited from $\mathbb{R}^{3}$,
meaning that $\tau_{\mathcal{M}}$ is also the same as the topology induced
by (the restriction of) $d$. Then, if a sequence $\{p_{n}\}$ on $\mathcal{M}$
is a Cauchy sequence with respect to $d_{\mathcal{M}}$, it is also a Cauchy
sequence with respect to $d$. Due to the completeness of $\mathbb{R}^{3}$
and the closedness of $\mathcal{M}$ in $\mathbb{R}^{3}$, there exists $\bar{p}\in \mathcal{M}$ such that $p_{n}\xrightarrow{d}\bar{p}$ (i.e., the sequence $\{p_n\}$ converges to $\bar{p}$ with respect to the metric $d$). Since $d_{\mathcal{M}}$ and $d$
induce the same topology on $\mathcal{M}$, this implies that $p_{n}\xrightarrow{d_{\mathcal{M}}}\bar{p}$ (i.e., the sequence $\{p_n\}$ converges to $\bar{p}$ with respect to the metric $d_{\mathcal{M}}$), ensuring the completeness of $(\mathcal{M},d_{\mathcal{M}})$.

The map $h:\mathcal{M}_{0}\rightarrow \mathcal{M}$ defined by 
\[
h(x,y,z)=\big(\sqrt{r(y)}\cdot x,\, y,\, \sqrt{r(y)}\cdot z\big)
\]
is a diffeomorphism between the pairs $(\mathcal{M}_{0},\mathcal{S}_{0})$ and $(\mathcal{M},\mathcal{S})$.
Here, we define $X$ to be the vector field on $\mathcal{M}$ related
to $X_{0}$ by $h$. That is, $X=h_{*}(X_{0})$, where $h_{*}: T \mathcal{M}_0 \to T \mathcal{M}$ is the tangent map. Let $\varphi$ be the flow of $X$ on $\mathcal{M}$. Then $h$ is
a conjugacy between the flows $\varphi_{0}$ and $\varphi$. That
is, the identity $h\circ\varphi_{0}=\varphi\circ h$ holds, or equivalently,
\begin{equation}
\varphi^{t}(p'')=h\circ\varphi_{0}^{t}\circ h^{-1}(p'')\label{eq:conjugacy of phi and phi_0}
\end{equation}
for all $p''\in \mathcal{M}$.

Note that for a point $p'=(x',y',z')$ on $\mathcal{M}_{0}$, the distance $d_{\mathcal{M}_{0}}(p',\mathcal{S}_{0})$ is \emph{exactly} the length of the minor arc on the circle $\mathcal{M}_{0}\cap\{(x,y,z) \in \mathbb{R}^3 : y=y'\}$
between $p'$ and $(0,y',1)$. Meanwhile, for a point $p''=(x'',y'',z'')=h(p')$
on $\mathcal{M}$, the distance $d_{\mathcal{M}}(p'',\mathcal{S})$ is \emph{no larger than} the length of the minor arc on the circle $\mathcal{M}\cap\{(x,y,z) \in \mathbb{R}^3 :  y=y''\}$ between $p''$ and $\big(0,y'',\sqrt{r(y'')}\big)$.
With $r(y)\leq1$, this implies 
\[
d_{\mathcal{M}}\big(h(p'),\mathcal{S}\big)\leq d_{\mathcal{M}_{0}}(p',\mathcal{S}_{0})
\]
for all $p'\in \mathcal{S}_{0}$. Combined with \eqref{eq:conjugacy of phi and phi_0},
it yields the following inequality:
\begin{equation}
d_{\mathcal{M}}\big(\varphi^{t}(p''),\mathcal{S}\big)=d_{\mathcal{M}}\big(h\circ\varphi_{0}^{t}\circ h^{-1}(p''),\mathcal{S}\big)\leq d_{\mathcal{M}_{0}}\big(\varphi_{0}^{t}\circ h^{-1}(p''),\mathcal{S}_{0}\big).\label{eq:d_M=000026d_M0}
\end{equation}
Since $h^{-1}$ maps $\tilde{\mathcal{D}}_{0} \defeq \{(x,y,z)\in \mathcal{M} : z>-1\}$ diffeomorphically to $\mathcal{D}_{0}$, it implies that as $t\rightarrow+\infty$,
$d_{\mathcal{M}}\big(\varphi^{t}(p),\mathcal{S}\big)\rightarrow0$ for all $p\in\tilde{\mathcal{D}}_{0}$.
However, if $\mathcal{S}$ is an attractor, then the domain of attraction of $\mathcal{S}$ should be
\[
\mathcal{D}=\tilde{\mathcal{D}}_{0}\cup\{\big(0,y,-\sqrt{r(y)}\big) : y>0\}.
\]
To see this, first note that for any point $p''=(x'',y'',z'')$ in $\{\big(0,y,-\sqrt{r(y)}\big) : y>0\}$, 
\[
\begin{aligned}\varphi^{t}(p'') &  & = & h\circ\varphi_{0}^{t}\circ h^{-1}(p'')\\
&  & = & h\circ\varphi_{0}^{t}(0,y'',-1)\\
&  & = & h(0,\gamma''(t),-1)\\
&  & = & \big(0,\gamma''(t),\sqrt{r\circ\gamma''(t)}\big),
\end{aligned}
\]
where $\frac{d\gamma''}{dt}>0$. Then from  (\ref{eq:dot_gamma>0}) we can deduce that $\gamma''(t)$ and $\frac{d\gamma''}{dt}$ both strictly increase with respect to $t$. Hence $d_{\mathcal{M}}(\varphi^{t}(p''),\mathcal{S})\leq\pi\sqrt{r\circ\gamma''(t)}\rightarrow0$
as $t\rightarrow+\infty$. Meanwhile, for any point $p\in \mathcal{M}-\mathcal{D}$,
i.e. $p=(0,y,-1)$ with $y\leq0$, $X|_{p}=h_{*}(X_{0} |_{p})=0$, and hence $p$ stays stationary under the flow $\varphi$. Therefore, if $p''\in \mathcal{M}$, then $\varphi^{t}(p)\xrightarrow{d_{\mathcal{M}}}\mathcal{S}$ as $t \to \infty$ if and only if $p''\in\mathcal{D}$. Since $\mathcal{D}$ contains circles in the
form $\{(x,y,z) \in \mathbb{R}^3 : x^{2}+z^{2}=r(y),y>0\}$ in $\mathcal{M}$, its fundamental group is
non-zero and hence is not homotopy equivalent to $\mathcal{S}$.

To show that this is a counterexample, it remains to prove that $\mathcal{S}$ is indeed a uniformly asymptotically stable manifold of
the system $(\mathcal{M},\varphi)$. Let 
\[
\mathcal{W} \defeq \{(x,y,z)\in \mathcal{M} : z>0\}\cup\{(x,y,z)\in \mathcal{M} : y>1\}.
\]
We will first show that $\mathcal{W}$ contains some $\alpha$-neighborhood
$\mathcal{N}_{\alpha}$ of $\mathcal{S}$ for some $\alpha>0$, and then show
that for each $\epsilon>0$, there exists some $T_{\epsilon}>0$ such
that $d_{\mathcal{M}}\big(\varphi^{[T_{\epsilon},+\infty)}(\mathcal{W}),\mathcal{S}\big)<\epsilon$. 

To see that $\mathcal{W}$ contains some $\alpha$-neighborhood of $\mathcal{S}$, we
only need to show that there is a positive distance between its complement
$\mathcal{W}^{c}$ and $\mathcal{S}$. Note that 
\[
\mathcal{W}^{c}=\{(x,y,z)\in \mathcal{M} : z\leq0,y\leq1\}=\mathcal{C}\cup \mathcal{K}
\]
with 
\[
\mathcal{C} \defeq \{(x,y,z)\in \mathcal{M} : z\leq0,y\leq-\pi\}
\]
and 
\[
\mathcal{K} \defeq \{(x,y,z)\in \mathcal{M} : z \leq 0,-\pi\leq y\leq1\}.
\]
Then $\mathcal{K}$ is compact and $\mathcal{C}$ is closed in $\mathcal{M}$, and $\mathcal{C}\cap \mathcal{S}$, $\mathcal{K}\cap \mathcal{S}$
are both empty. Since $(\mathcal{M},g_{\mathcal{M}})$ is a complete Riemannian manifold
with the distance $d_{\mathcal{M}}$, it holds that $d_{\mathcal{M}}(\mathcal{S},\mathcal{K})>0$ as a consequence of the disjointedness of a closed subset and a compact subset. To see that $d_{\mathcal{M}}(\mathcal{S},\mathcal{C})>0$, note that $d_{\mathcal{M}}(\mathcal{S}_{y\leq0},\mathcal{C}) = \pi / 2$ and
$d_{\mathcal{M}}(\mathcal{S}_{y\geq0},\mathcal{C})\geq\pi$, where $\mathcal{S}_{y\leq0} \defeq \mathcal{S}\cap\{ (x,y,z) \in \mathbb{R}^3 : y\leq0\}$ and
$\mathcal{S}_{y\geq0} \defeq \mathcal{S}\cap\{ (x,y,z) \in \mathbb{R}^3 : y\geq0\}$. Then for any $0<\alpha<\min\{d_{\mathcal{M}}(\mathcal{S},\mathcal{K}),d_{\mathcal{M}}(\mathcal{S},\mathcal{C})\}$, there holds $\mathcal{N}_{\alpha}\subset \mathcal{W}$.

Now we proceed to show that for any $\epsilon>0$, there exists $T_{\epsilon}>0$. Denote by $\mathcal{W}_{z>0}$ the set $\{(x,y,z)\in \mathcal{M}: z>0\}$ and
by $\mathcal{W}_{y>1}$ the set $\{(x,y,z)\in \mathcal{M} : y>1\}$. Then $\mathcal{W}=\mathcal{W}_{z>0}\cup \mathcal{W}_{y>1}$.
Note that for each point $p=(x,y,z)\in \mathcal{W}_{y>1}$, $X |_{p}$
takes the form $a_{p}\frac{\partial}{\partial x}+e^{-\frac{1}{y}}\frac{\partial}{\partial y}+c_{p}\frac{\partial}{\partial z}$,
and therefore, $\mathcal{W}_{y>1}$ is an invariant open set of the system $(\mathcal{M},\varphi)$.
It holds that $\varphi^{t}(p)=\big(x_{t},y_{t},z_{t}\big)$ with $\frac{dy_{t}}{dt}>e^{-1}$
for any $p\in \mathcal{W}_{y>1}$. Choose $T''$ to be some positive number
large enough such that $r(e^{-1}\cdot T'')<(\epsilon / \pi)^{2}$.
Then for any $t\geq T''$ and $p\in \mathcal{W}_{y>1}$, it holds that $r(y_{t})<r(e^{-1}\cdot T'')$
and therefore $d_{\mathcal{M}}\big(\varphi^{t}(p),\mathcal{S}\big)\leq\pi\sqrt{r(y_{t})}<\epsilon$.
To see that the points in $\mathcal{W}_{z>0}$ converge uniformly towards $\mathcal{S}$,
first note that $h^{-1}(\mathcal{W}_{z>0})=\mathcal{W}'_{z>0}=\mathcal{M}_{0}\cap\{(x,y,z) \in \mathbb{R}^3 : z>0\}$. Then
combined with \eqref{eq:d_M=000026d_M0}, there holds $d_{\mathcal{M}}\big(\varphi^{[T',+\infty)}(\mathcal{W}_{z>0}),\mathcal{S}\big)<\epsilon$.
Finally, one only needs to choose $T_{\epsilon}$ to be $\min\{T',T''\}$ and the whole argument is complete.

\section{Conclusion} \label{sec_conclu}
In this paper, we have revisited Wilson's theorem (i.e., Theorem 3.4 in \cite{wilson1967structure}) about the relation between the domain of attraction of an attractor and its tubular neighborhood. Specifically, we show with detailed and rigorous proofs that the domain of attraction of a \emph{compact} asymptotically submanifold of a finite-dimensional smooth manifold of a continuous dynamical system is homeomorphic to its tubular neighborhood. We emphasize that the compactness of the attractor is crucial, without which Wilson's theorem cannot hold. This is shown by two counterexamples where the attractor is not compact and the state space is either complete or incomplete. 

\bibliographystyle{plain}
\bibliography{refnote}

\end{document}